\documentclass[10pt,a4paper,twoside,leqno]{article}
\usepackage{amsfonts,amssymb}

\usepackage{amscd}

\vspace{0.5cm}

 4
\font\ebf=cmbx8
\font\erm=cmr8

\usepackage{graphicx}

\begin{document}
\begin{center}
	\noindent { \textsc{ Lucky 13-th Exercises on Stirling-like numbers and Dobinski-like formulas}}  \\ 
	\vspace{0.3cm}
	\vspace{0.3cm}
	\noindent Andrzej Krzysztof Kwa\'sniewski \\
	\vspace{0.2cm}
	\noindent {\erm Member of the Institute of Combinatorics and its Applications  }\\
{\erm High School of Mathematics and Applied Informatics} \\
	{\erm  Kamienna 17, PL-15-021 Bia\l ystok, Poland }\noindent\\
	\noindent {\erm e-mail: kwandr@gmail.com}\\
	\vspace{0.4cm}
\end{center}

\vspace{0.1cm}

\noindent {\ebf Summary}

\noindent {\small Extensions of the  Stirling numbers of the second kind
and Dobinski-like formulas are proposed in a series of exercises
for graduetes. Some of  these new formulas recently discovered  by me are
to be found in A.K.Kwa\'sniewski's source paper [1]. These extensions naturally
encompass the well known $q$-extensions.The indicatory references are to point at
a part of the vast domain of the foundations of computer science in ArXiv affiliation
 noted as CO.cs.DM.

\vspace{0.2cm}

\noindent MCS numbers:  05A40, 11B73, 81S99
\vspace{0.2cm}

\noindent Keywords: {umbral calculus, extended Stirling numbers, Dobinski type identities}

\vspace{0.2cm}

\noindent  affiliated to The Internet Gian-Carlo Polish Seminar:

\noindent \emph{http://ii.uwb.edu.pl/akk/sem/sem\_rota.htm}

\vspace{0.1cm}

\noindent Published in:    Proc. Jangjeon Math. Soc.  Vol. \textbf{11}  No 2, (\textbf{2008}),137-144 
\vspace{0.2cm}

\subsection*{1. In the $q$-extensions realm}

\textbf{Ex.1} \textbf{ Recall and prove} it again occasionally noting that The Number $44$ is a magic number in Polish 
Poetry and this had had an implementation in quite recent Polish history in $1968$ students revolutionary riots
for independence and freedom of thinking under communist regime. Well,  then forty four years
ago Gian-Carlo Rota [4] proved that the exponential
generating function for Bell numbers $B_n$ is of the form

\begin{equation}\label{eq1}
     \sum_{n=0}^\infty \frac {x^n}{n!}B_n = \exp(e^x-1)
\end{equation}
 using the linear functional  \textit{L } such that
\begin{equation}\label{eq2}
     L(X^{\underline{n}})=1  , \qquad  n\geq 0
\end{equation}

\noindent Then Bell numbers (see: formula (4)  in [4])  are defined by
\begin{equation}\label{eq3}
          L(X^n)=B_n  ,\qquad   n \geq 0
\end{equation}
The above formula is exactly the Dobinski formula [5] if $L$
is interpreted as the average value functional for the random
variable $X$  with the Poisson distribution with $L(X) = 1$.
Recall it and prove it again.

\vspace{1mm}

\noindent \textbf{Ex.2}   \textbf{Recalland prove}:

\noindent The two standard [12],see also [13-17] $q$-extensions Stirling numbers of the
second kind  are defined  by

\begin{equation}\label{eq4}
x_q^n=\sum_{k=0}^{n}\Big\{{n \atop k}\Big\}_q  x_q^{\underline k},
\end{equation}

\noindent where $ x_q=\frac{1-q^x}{1-q}$ and $x_q^{\underline
k}=x_q(x-1)_q...(x-k+1)_q $ , which corresponds to the $\psi$
sequence choice in the $q$-Gauss form
$\langle\frac{1}{n_q!}\rangle_{n\geq 0}$ and $q^\sim$-Stirling
numbers

\begin{equation}\label{eq5}
x^n=\sum_{k=0}^{n}\Big\{{n \atop k}\Big\}^\sim_q \chi_{\underline
k}(x)
\end{equation}

\noindent where $\chi_{\underline k}(x)= x(x-1_q)(x-2_q)...(x-[k-1]_q)$

\noindent Note that these formulae become the  usual unextended Stirling
numbers of the second kind formulae when the subscript $q$ is
removed.

\vspace{1mm}

\noindent \textbf{Ex.3}  \textbf{Recall and prove} it again:

\vspace{1mm}

\noindent For these two classical by now $q$-extensions of  Stirling numbers
of the second kind  - the "$q$-standard" recurrences hold
respectively:

$$ \Big\{{{n+1}\atop k}\Big\}_q =
\sum_{l=0}^{n}{n \choose l}_q q^l\Big\{{l\atop {k-1}}\Big\}_q ;
n\geq 0 , k\geq 1,$$

$$ \Big\{{{n+1}\atop k}\Big\}^\sim_q =
\sum_{l=0}^{n}{n \choose l}_q q^{l-k+1}\Big\{{l\atop
{k-1}}\Big\}^\sim_q  ;  n\geq 0 , k\geq 1.$$

\noindent Show that these formulae become the  usual unextended  Stirling
numbers of the second kind formulae when the subscript $q$ is
removed  and the number $q$  is put equal to one.

\vspace{1mm}

\noindent \textbf{Ex.4}  \textbf{Recall and prove} it again:

\vspace{1mm}

\noindent From the above it follows immediately that corresponding
$q$-extensions of  $B_n$ Bell numbers satisfy respective
recurrences:

$$ B_q(n+1) =
\sum_{l=0}^{n}{n \choose l}_q q^l B_q(l) ; n\geq 0 ,$$

$$ B^\sim_q(n+1) =
\sum_{l=0}^{n}{n \choose l}_q q^{l-k+1} \overline {B} ^\sim_q(l)  ;
n\geq 0 $$ where
$$ \overline {B}^\sim_q(l)=
\sum_{k=0}^{l} q^k\Big\{{l\atop k}\Big\}^\sim_q .
$$
Show that these formulae become the  usual unextended  Stirling
numbers of the second kind formulae when the subscript $q$ is
removed  and the number $q$  is put equal to one.

\vspace{1mm}

\noindent \textbf{Ex.5}  \textbf{Recall and prove} it again:

\vspace{1mm}

\noindent Recursions for both \textit{inversion} $q$-Bell numbers and
\textit{inversion} $q$-Stirling  numbers of the second kind are
not difficult to be derived. Also in a natural way the
\textit{inversion} $q$-Stirling numbers of the second kind from
[16] satisfy a $q$-analogue of the standard recursion for Stirling
numbers of the second kind to be written via mnemonic adding "$q$"
subscript to the binomial and second kind Stirling symbols in the
the standard recursion formula  i.e.

$$ \Big\{{{n+1}\atop k}\Big\}^{inv}_q =
\sum_{l=0}^{n}{n \choose l}_q \Big\{{{n-l}\atop
{k-1}}\Big\}^{inv}_q  ;  n\geq 0 , k\geq 1.$$
Another $q$-extended Stirling numbers much different from Carlitz
"$q$-ones" were introduced in the reference [19],see [1,20,28].

\noindent The \textit{cigl}-$q$-Stirling numbers of the second kind  are
expressed in terms of $q$-binomial coefficients and $q =1$
Stirling numbers of the second kind [16,17],( see [1] for more references) as follows

$$ \Big\{{{n+1}\atop k}\Big\}^{cigl}_q =
\sum_{l=0}^{n}{n \choose l}_q q^{{n-l+1 \choose 2}}\Big\{{{n-l}\atop
{k-1}}\Big\}^{cigl}_q ; n\geq 0 , k\geq 1.$$
The corresponding \textit{cigl}-$q$-Bell numbers recently have
been equivalently defined via \textit{cigl}-$q$-Dobinski formula
[20,28] - which now in more adequate notation reads :

$$L(X^{\overline{q^n}})=\overline{B}_n(q),\qquad n\geq 0,
X^{\overline{q^n}}\equiv X(X+q-1)...(X-1+q^{n-1}).$$
The above \textit{cigl}-$q$-Dobinski formula is interpreted as the average
of this specific  $n-th$ \textit{cigl}-$q$-power random variable
$X^{\overline{q^n}} $ with the $q = 1$ Poisson distribution such
that $ L(X)=1 $.
\vspace{0.1cm}

\noindent For that to see use the identity by Cigler [19]
$$
x(x-1+q)...(x-1+q^{n-1})=\sum_{k=0}^{n}\Big\{ {n \atop
k}\Big\}^{cigl}_q  x^{\underline k}.
$$

\noindent Note that these formulae become the  usual unextended  Stirling
numbers of the second kind formulae when the subscript $q$ is
removed and the number $q$  is put equal to one.

\vspace{1mm}

\subsection*{2. Beyond the $q$-extensions realm}

The further consecutive umbral extension of Carlitz-Gould
$q$-Stirling numbers   $\Big\{{n\atop k}\Big\}_q $  and
$\Big\{{n\atop k}\Big\}^\sim_q $   is realized two-fold way - one
of which leads to a surprise (?) in  contrary perhaps to the other way. \\

\textbf{The first} "easy way" consists in almost mnemonic
sometimes replacement of $q$ subscript by $\psi$ after having
realized that via equation (5) we are dealing with the specific
case of Comtet numbers [1] i.e. now we have

\begin{equation}\label{eq6}
x^n=\sum_{k=0}^{n}\Big\{{n \atop k}\Big\}^\sim_{\psi}
\psi_{\underline k}(x)
\end{equation}

\noindent where $\psi_{\underline k}(x)=
x(x-1_{\psi})(x-2_{\psi})...(x-[k-1]_{\psi}).$ Show that these
formulae become the  usual unextended  Stirling numbers of the
second kind formulae when the subscript $\psi$ is removed.

\vspace{1mm}

\noindent As a consequence we have "for granted" the following: \vspace{2mm}

\vspace{1mm}

\noindent \textbf{Ex.6}  \textbf{Recall standard and prove} its extension:

\vspace{1mm}

\begin{equation}\label{eq7}
\Big\{{{n+1}\atop k}\Big\}^\sim_{\psi} = \Big\{{n\atop
{k-1}}\Big\}^\sim_{\psi} + k_{\psi}\Big\{{n\atop
k}\Big\}^\sim_{\psi} ;\quad n\geq 0 , k\geq 1 ;
\end{equation}

\noindent where \quad $\Big\{{n\atop 0}\Big\}^\sim_{\psi}=
\delta_{n,0},\quad\Big\{{n\atop k}\Big\}^\sim_{\psi}=0 ,\quad k>n
;\quad $\quad and the recurrence for ordinary generating function
reads

\vspace{1mm}

\begin{equation}\label{eq8}
G^\sim_{k_{\psi}}(x)=\frac{x}{1-k_{\psi}}G^\sim_{k_{\psi}-1}(x) ,
\quad k\geq 1
\end{equation}

\noindent where naturally
$$G^\sim_{k_{\psi}}(x)= \sum_{n\geq 0}\Big\{{n\atop
k}\Big\}^\sim_{\psi}x^n     ,\quad k\geq 1 $$

\vspace{1mm}

\noindent from where one infers \vspace{5mm}

\begin{equation}\label{eq9}
G^\sim_{k_{\psi}}(x)=\frac{x^k}{(1-1_{\psi}x)(1-2_{\psi}x)...(1-k_{\psi}x)}
\quad, \quad k\geq 0
\end{equation}

\vspace{1mm}

\noindent hence we arrive in the standard extended  text-book way [21]
at the following explicit formula

\begin{equation}\label{eq10}
\Big\{{n\atop k}\Big\}^\sim_{\psi} = \frac
{r_{\psi}^n}{k_{\psi}!}\sum_{r=1}^{k}(-1)^{k-r}{k_{\psi} \choose r_{\psi}};
\quad n,k\geq 0 .
\end{equation}
\noindent Show that these formulae become the  usual unextended  Stirling
numbers of the second kind formulae when the subscript $\psi$ is
removed.

\vspace{1mm}

\noindent Expanding the right hand side of the corresponding equation above
results in another explicit formula for these $\psi$-case Comtet
numbers [1] i.e. we have

\vspace{1mm}

\noindent \textbf{Ex.7}  \textbf{Recall standard and prove} its extension:

\begin{equation}\label{eq11}
\Big\{{n\atop k}\Big\}^\sim_{\psi} = \sum_{1\leq i_1 \leq
i_2\leq...\leq i_{n-k}\leq
k}(i_1)_{\psi}(i_2)_{\psi}...(i_{n-k})_{\psi}; \quad n,k\geq 0 .
\end{equation}

\vspace{1mm}

\noindent or equivalently  (compare with [13], see [12,14,15])

\vspace{1mm}

\begin{equation}\label{eq12}
\Big\{{n\atop k}\Big\}^\sim_{\psi}=\sum_{d_1+ d_2+...+d_k =
n-k,\quad d_i\geq
0}1_{\psi}^{d_1}2_{\psi}^{d_2}...k_{\psi}^{d_k};\quad n,k\geq 0 .
\end{equation}

\vspace{1mm}

\noindent $\psi^\sim$-\textbf{Stirling numbers} of the second kind being
defined equivalently by (10) , (14), (15)  or (16) yield
$\psi^\sim$-\textbf{Bell numbers}
$$ B^\sim_n(\psi)=\sum_{k=0}^{n} \Big\{{n\atop
k}\Big\}^\sim_{\psi} ,\qquad n\geq 0 .$$

\noindent Show that these formulae become the  usual unextended  Stirling
numbers of the second kind formulae when the subscript $\psi$ is
removed.

\vspace{1mm}

\noindent \textbf{Ex.8}   \textbf{Recall standard and prove} its extension:

\noindent Adapting the standard text-book method  [21] we have for two
variable ordinary generating function for $\Big\{{n\atop
k}\Big\}^\sim_{\psi}$ Stirling numbers of the second kind and the
$\psi$-exponential generating function for $ B^\sim_n(\psi)$  Bell
numbers the following  formulae

\begin{equation}\label{eq13}
C^\sim_{\psi}(x,y) = \sum_{n\geq 0} A^\sim_n (\psi,y)x^n ,
\end{equation}
where the $\psi$- exponential-like polynomials $ A^\sim_n
(\psi,y)$
$$ A^\sim_n (\psi,y)=\sum_{k=0}^{n} \Big\{{n\atop
k}\Big\}^\sim_{\psi}y^k$$ do satisfy the recurrence

$$ A^\sim_n (\psi,y)=  [y(1+\partial_{\psi}]A^\sim_{n-1}(\psi,y) \qquad n\geq 1 ,$$

\noindent hence
$$ A^\sim_n (\psi,y)=  [y(1+\partial_{\psi}]^n 1,\quad \qquad n\geq 0 ,$$
where  the linear operator $\partial_{\psi}$ acting on the algebra
of formal power series is being called (see: [1,2,3,24,25,31])
the "$\psi$-derivative" and $\partial_{\psi} y^n =
n_{\psi}y^{n-1}.$

\noindent Show that these formulae become the  usual unextended Stirling
numbers of the second kind formulae when the subscript $\psi$ is
removed.

\vspace{1mm}

\noindent \textbf{Ex.9} \textbf{Recall standard and prove} its extension:

\noindent The $\psi$-exponential generating function $ B^{\sim}_{\psi}(x)=
\sum_{n\geq 0}B^\sim_n(\psi)\frac{x^n}{n_{\psi}!} $

\noindent for $ B^{\sim}_n(\psi)$ Bell numbers - after cautious adaptation
of the method from the Wilf`s generatingfunctionology book  [21]
can be seen to be  given by the following   formula

\begin{equation}\label{eq14}
B^{\sim}_{\psi}(x)= \sum_{r\geq 0}\frac{1}{\epsilon(\psi,r)}
\frac{e_{\psi}[r_{\psi}x]}{r_{\psi}!}
\end{equation}

\noindent where (see: [1,2,3,24,25,31]) 
$$e_{\psi}(x) =
\sum_{n\geq 0}\frac{x^n}{n_{\psi}!}$$

\noindent while
\begin{equation}\label{eq15}
\epsilon(\psi,r)=\sum_{k=r}^{\infty} \frac{(-1)^{k-r}}{(k_{\psi}-
r_{\psi})!}
\end{equation}
\noindent and for the $\psi$-extensions the Dobinski like formula here now
reads

\begin{equation}\label{eq16}
B^{\sim}_n (\psi)= \sum_{r\geq 0}\frac{1}{\epsilon(\psi,r)}
\frac{r_{\psi}^n}{r_{\psi}!}.
\end{equation}
\noindent Show that these formulae become the  usual unextended  Stirling
numbers of the second kind formulae when the subscript $\psi$ is
removed.

\vspace{1mm}

\noindent \textbf{Ex.10}  \textbf{Recall standard and prove} its extension:

\noindent In the case of  Gauss $q$-extended  choice of
$\langle\frac{1}{n_q!}\rangle_{n\geq 0}$ admissible sequence of
extended umbral operator calculus  equations  (19) and (20) take
the form

\begin{equation}\label{eq17}
\epsilon(q,r)=\sum_{k=r}^{\infty}
\frac{(-1)^{k-r}}{(k-r)_q!}q^{- {r \choose 2}}
\end{equation}
\noindent and the $q^\sim$-Dobinski  formula is given by

\begin{equation}\label{eq18}
B^{\sim}_n (q)= \sum_{r\geq 0}\frac{1}{\epsilon(q,r)}
\frac{r_q^n}{r_q!},
\end{equation}
\noindent which for $ q=1$ becomes the Dobinski formula from 1887 [5].

\vspace{1mm}

\noindent \textbf{Ex.11} \textbf{Recall standard and prove} its extension:

\noindent In a dual inverse way  we define the $\psi^\sim$ Stirling numbers
of the first kind as coefficients in the following expansion

\begin{equation}\label{eq19}
\psi_{\underline k}(x)=\sum_{r=0}^{k}\left[ {{\begin{array}{*{20}c} {k}\\
{r}\end{array}} } \right]^\sim_\psi x^r
\end{equation}

\noindent where -  recall $\psi_{\underline k}(x)=
x(x-1_{\psi})(x-2_{\psi})...(x-[k-1]_{\psi});\quad$( attention:
see equations (10)-(16)\quad in [7,8] and note the difference with
the present definition).

\noindent Therefore from the above we  infer that

\begin{equation}\label{eq20}
\sum_{r=0}^{k}\left[ {{\begin{array}{*{20}c} {k}\\
{r}\end{array}} } \right]^\sim_\psi \Big\{{r\atop
l}\Big\}^\sim_{\psi}= \delta_{k,l}.
\end{equation}
\noindent Show that these formulae become the  usual unextended  Stirling
numbers of the second kind formulae when the subscript $\psi$ is
removed.

\vspace{1mm}

\noindent \textbf{Ex.12} \textbf{Recall standard and prove} its extension:

\noindent Consider now another Stirling-like numbers (as expected  Whitney numbers [31,1,32]) which are 
natural counterpart to $\psi^\sim$-Stirling numbers of the second kind. These  are $\psi^c$-
Stirling numbers  of the first kind defined here down  as coefficients in the following
expansion (upperscript "\textit{c}" is used because of cycles in non-extended case).

\begin{equation}\label{eq21}
\psi_{\overline k}(x)=\sum_{r=0}^{k}\left[ {{\begin{array}{*{20}c} {k}\\
{r}\end{array}} } \right]^c_\psi x^r
\end{equation}

\noindent where -  now  $\psi_{\overline k}(x)=
x(x+1_{\psi})(x+2_{\psi})...(x+[k-1]_{\psi});\quad$

\noindent Show that these formulae become the  usual unextended  Stirling
numbers of the second kind formulae when the subscript $\psi$ is
removed.

\vspace{1mm}

\noindent \textbf{Ex.13} \textbf{Recall standard and prove} its extension:

\noindent Show that the definition here upstairs above of the $q$-Stirling
numbers of the second kind $\Big\{{n \atop k}\Big\}_q$  is
equivalent with the definition by recursion

\begin{equation}\label{eq22}
\Big\{{{n+1}\atop k}\Big\}_q = q^{k-1}\Big\{{n\atop {k-1}}\Big\}_q
+ k_q\Big\{{n\atop k}\Big\}_q ;\quad n\geq 0 , k\geq 1 ;
\end{equation}

\noindent where \quad $\Big\{{n\atop 0}\Big\}_q =
\delta_{n,0},\quad\Big\{{n\atop k}\Big\}_q=0 ,\quad k>n \quad $

\vspace{1mm}

\noindent Show that these in turn \quad (just use the $Q$-Leibnitz rule
[2,3,24,25,31] for Jackson derivative $\partial_q$) \quad  are
equivalent to

\begin{equation}\label{eq23}
( \hat{x} \partial_q)^n =\sum_{k=0}^{n}\Big\{{n \atop k}\Big\}_q
\hat{x}^k
\partial_q^k
\end{equation}

\noindent where \quad $\Big\{{n\atop 0}\Big\}_q =
\delta_{n,0},\quad\Big\{{n\atop k}\Big\}_q=0 ,\quad k>n \quad .$

\noindent Here $\hat{x}$ denotes the multiplication by the argument of a
function.

\vspace{1mm}

\noindent Show that these formulae become the  usual unextended
Stirling numbers of the second kind formulae when the subscript
$q$ is removed.

\vspace{3mm}
\noindent  Consult [33-35] for some new open problems arising in related the domain of the so called 
cobewb posets and their acyclic digraphs representatives which had served the present author to discover
a joint combinatorial interpretation for all $F$-nomial coefficients. These family encompasses binomial and $q$-Gaussian coefficient,
 Fibonomial coefficient, Stirling numbers  of both kinds and all classical $F-nomial$ coefficients hence specifically incidence
coefficients of reduced incidence algebras of full binomial type and Whitney numbers are given the joint cobweb combinatorial interpretation also
[36,37].

\vspace{3mm}

\noindent \textbf{Acknowledgements} The author appraises much Maciej Dziemia\'nczuk Gda\'nsk Univerity Student's assistance including $TeX-nology$ aid.

The author expresses his gratitude  also   Dr Ewa Krot-Sieniawska for her several years' cooperation and vivid application  of the alike
material deserving  Students' admiration for her being such a comprehensible and reliable  Teacher before \textbf{she was fired  by Bialystok University authorities exactly on the day she had defended  Rota and cobweb posets related dissertation with distinction}.

\end{document}